\documentclass[12pt]{article}
\usepackage[all]{xy}
\usepackage{wrapfig}

\usepackage{epic}
\usepackage{amsmath}[1996/11/01]
\usepackage{amssymb,amsthm,amsfonts,latexsym,epsfig,graphics}
 \def\beql#1#2\eeql{\begin{equation}\label{#1}#2\end{equation}}

\advance \textheight by -1 \topmargin
\advance \textheight by -1 \headheight
\advance \textheight by -1 \headsep
\advance \textheight by -2 \footskip
\vsize = \textheight
\hsize = \textwidth
\newcommand{\meet}{\cap }

\DeclareMathOperator{\WQ}{WQ}

\DeclareMathOperator{\rk}{rk}

\DeclareMathOperator{\soc}{soc}

\DeclareMathOperator{\End}{End}

\DeclareMathOperator{\Char}{char}

\DeclareMathOperator{\GL}{GL}

\DeclareMathOperator{\Hom}{Hom}

\DeclareMathOperator{\Sp}{Sp}
\DeclareMathOperator{\id}{id}

\theoremstyle{plain}
\newtheorem{theorem}{Theorem}
\newtheorem{lemma}[theorem]{Lemma}

\newtheorem{proposition}[theorem]{Proposition}
\newtheorem{corollary}[theorem]{Corollary}
\newtheorem{definition}[theorem]{Definition}
\newtheorem{remark}[theorem]{Remark}
\theoremstyle{remark}

\numberwithin{theorem}{section}

\newcommand{\disj}{\stackrel{.}{\cup}}

\newcommand{\F}{{\mathbb{F}}}

\newcommand{\A}{{\mathsf{A}}}

\newcommand{\C}{{\mathsf{C}}}
\newcommand{\D}{{\mathsf{D}}}

\renewcommand{\em}{\sf}

\title{Equivariant quadratic forms in characteristic 2}
\author{Gabriele Nebe and Richard Parker} 
\date{}
 
\begin{document}
\maketitle

 {\sc Abstract.} 
 Let $G$ be a finite group and $K$ a finite field of characteristic $2$. 
 Denote by $t$ the $2$-rank of the commutator factor group $G/G'$ and 
 by $s$ the number of isomorphism classes of self-dual simple $KG$-modules. 
 Then the 
 Witt group of equivariant quadratic forms $\WQ (K,G)$ is isomorphic
 to an elementary abelian $2$-group of rank $s+t$.

\section{Introduction}
\label{sec:intro}

Witt groups of quadratic and Hermitian forms have intensively 
been studied by various authors. 
 In particular the paper 
 \cite{Dress} lies the foundations for a theory of 
 equivariant quadratic forms for finite groups. 
 Most approaches in the literature deal with bilinear 
 or Hermitian forms.
 The textbook \cite{Rhiem} investigates orthogonal 
 representations and equivariant Witt groups for fields of characteristic not 2 
 (see also Section \ref{odd} for an explicit description). 
Equivariant quadratic 
forms over fields of characteristic 2 require adapted methods 
 as developed for instance in \cite{SinWillems}. 
 The aim of the present paper is to give an
 explicit description of the equivariant Witt groups $\WQ(K,G)$ 
 of quadratic forms for finite groups $G$ and finite fields $K$ 
 of characteristic 2. 
 The group elements of $\WQ(K,G)$ are equivalence classes 
 $[(V,Q)]$ of quadratic $KG$-modules $(V,Q)$. 
 Here $Q : V\to K$ denotes a $G$-invariant quadratic form on
 the $KG$-module $V$ that is non-degenerate, 
 i.e. the radical of its polarization (see equation \eqref{polar} below) 
 is zero. 
 Addition in $\WQ(K,G)$ is defined via the orthogonal direct
 sum of representatives. 
 By Theorem \ref{anisoWitt} each class in the Witt group has 
 a unique anisotropic representative, i.e. an equivariant 
 quadratic form $(V,Q)$ for which the restriction of $Q$ to
 any non-zero submodule is non-zero. 
 The main result of the present note is the following theorem. 

 \begin{theorem}
Let 
$K$ be  a finite field of characteristic $2$ and $G$ be a finite group. 
 Let $s$ denote the number of isomorphism classes of self-dual 
 simple	$KG$-modules (including the trivial module) 
 and let $t$ denote the $2$-rank of $G/G'$.
	 Then the equivariant Witt group 
	 $\WQ(K,G)$ is isomorphic to $C_2^{s+t}$,
  the elementary abelian $2$-group of rank $s+t$.
 \end{theorem}

 Generators of $\WQ(K,G)$ can be constructed as anisotropic 
 equivariant quadratic forms:

 For a simple $KG$-module $V$ admitting a 
 non-degenerate  $G$-invariant 
 quadratic form $Q$, this form is unique up to $G$-isometry 
 (see Theorem \ref{simple} (e)) and 
 defines the generator $[(V,Q)]$ 
 of the Witt group. 

 The trivial $KG$-module $T$ has dimension 1, is self-dual,
 but does not carry a 
 non-degenerate quadratic form. To this module we associate 
 the generator $[N(K)]$ where $N(K)$ 
 is the unique $2$-dimensional anisotropic quadratic space over $K$ 
 with trivial $G$-action. 

 The other simple self-dual $KG$-modules $W$ 
 carry a unique non-degenerate symplectic $G$-equivariant form 
 but no equivariant non-degenerate quadratic form. 
This yields a group homomorphism $G\to \Sp(W)$ into the symplectic
group on $W$. 
 Using the isomorphism $\Sp_{2m} (K) \cong O_{2m+1} (K) $ 
 and an embedding of the $2m+1$-dimensional semi-regular quadratic $K$-space
 into a non-degenerate quadratic space of dimension $2m+2$ (and maximal Witt index) 
we associate to $W$ the {\em quadratic envelope} of $W$ of type $+$, 
$[R^+(W)] $ (see Definition \ref{orthogonalisation})
as a generator of $\WQ(K,G)$. 

 The last set of generators is defined by the epimorphisms $\tau :G\to C_2$: 
 The quadratic $KG$-module $R^+(\tau ) $ has a basis $(b_1,b_2)$ which 
 is permuted by $\tau (G)$ and quadratic form $Q(a_1b_1+a_2b_2) = a_1a_2$. 
 For a basis $(\tau _1,\ldots , \tau _t) $ of $\Hom (G,C_2)$
 we get the $t$ classes $[R^+(\tau _j)] $ as additional generators. 

 For the proof we define group homomorphisms $\A $ and $\C $ 
 on $\WQ(K,G)$ in Section \ref{ACD}, 
 where $\A $ maps $[(V,Q)]\in \WQ(K,G)$
 to its class in the Witt group $\WQ(K)$ of quadratic $K$-spaces 
 and $\C $ takes those composition 
 factors occurring in $V$  with odd multiplicity. 
 Then the intersection of the kernels of $\A $  and $\C $ is 
 generated by the quadratic forms $[R^+(\tau _j)]$ above and hence
 isomorphic to $\Hom (G,C_2)$. 

 This research is funded under
Project-ID 286237555 -- TRR 195 -- by the
Deutsche Forschungsgemeinschaft (DFG, German Research Foundation).
We thank the referees for their helpful reports. 

\section{Witt groups of equivariant quadratic forms.} 

This section recalls the definition of Witt groups 
and exposes the short general argument that 
every class in the Witt group of equivariant quadratic forms 
over fields has a unique anisotropic representative. 
The results of this section are well known and can be found in
many textbooks, see for instance \cite[Chapter 7]{Scharlau}, \cite[Kapitel III]{Kneser}, or \cite{Dress}. 

\subsection{The Witt group of equivariant quadratic forms} 

Let $G$ be a finite group and $K$ be an arbitrary field. 
 An {\em equivariant quadratic form} $(V,Q)$ for $G$ consists of a right $KG$-module 
 $V$ together with a non-degenerate $G$-invariant quadratic form $Q:V\to K$. 
 Then the {\em polarization} $B_Q$ of $Q$ is defined as 
 the  $G$-invariant symmetric bilinear form given by
 \begin{equation}\label{polar}
B_Q(v,w) = Q(v+w)  - Q(v) - Q(w) \mbox{ for all } v,w\in V .
 \end{equation}
 The condition that $Q$ is non-degenerate is defined via the non-degeneracy
 of its polarization, the radical $V^{\perp }$ of $B_Q$ is $\{ 0 \}$.
 A submodule 
 $U\leq V$ is called {\em isotropic} if 
 $Q(U) = \{0 \}$. 

Any non-degenerate $G$-invariant bilinear form $B:V\times V \to K$ 
 yields a $KG$-isomorphism between $V$ and its 
 {\em dual module} $V^{\vee} = \Hom_{K}(V,K)$, in particular 
$V$ is a self-dual $KG$-module. 
 For a $KG$-submodule $U$ of $(V,B)$ the orthogonal space 
 $$U^{\perp }:= \{ v\in V \mid B(v,u) = 0 \mbox{ for all } u\in U \} $$
 is again a $KG$-submodule of $V$ and $V/U^{\perp} \cong U^{\vee }$. 

  \begin{definition}  \label{Witteq}
	  An equivariant quadratic form $(V,Q)$ is called 
	  {\em metabolic} if there is 
   an isotropic submodule 
  $U\leq V$ with $U=U^{\perp }$. 
Two quadratic forms $(V,Q)$ and $(W,Q')$ are called {\em Witt equivalent},
 if the orthogonal sum $(V,Q) \perp (W,-Q')$ is metabolic. 
 The equivalence classes $[(V,Q)]$ of non-degenerate equivariant 
 quadratic forms form an abelian group with orthogonal sum
 as addition, called the 
	  {\em Witt group} $\WQ(K,G)$
  of equivariant quadratic forms for $KG$.
  An equivariant quadratic form $(V,Q)$ is called 
  {\em anisotropic} if it does not contain a non-zero isotropic submodule. 
  \end{definition}

As a referee pointed out it is more common to define 
Witt equivalence by $(V,Q) \sim (W,Q')$ if and only if there
are metabolic modules $(N,H)$, $(M,H')$ such that 
$(V,Q) \perp (N,H) \cong (W,Q') \perp (M,H')$. 
In our situation this notion of Witt equivalence and the one 
in Definition \ref{Witteq} are equivalent. This is shown for 
Hermitian torsion $\pi $-spaces in \cite[Lemma 4.2]{Dress} and 
for Hermitian forms over finite algebras in \cite[Section 4.1]{Annika}. 
The proof can be taken almost literally also for equivariant quadratic forms. 
For convenience of the reader we include a proof that is based on the
following two lemmas:

\begin{lemma} \label{incl}
	Assume that $(V,Q)$ is metabolic and let $U\leq V$ be a $KG$-submodule 
	of $V$ such that 
	$Q(U) = \{ 0 \}$. Then there is a maximal isotropic $KG$-submodule
	$M = M^{\perp} $ of $V$ that contains $U$.
\end{lemma} 

\begin{proof} (see \cite[Lemma 4.1.4]{Annika})
	Let $N= N^{\perp} \leq V$ be a maximal isotropic $KG$-submodule. 
	Put 
	$$M := (N\cap U^{\perp}) + U .$$ 
	Then $Q(M) = \{0 \}$ and 
	$$M^{\perp} = ((N\cap U^{\perp}) + U  ) ^{\perp} = 
	(N+U) \cap U^{\perp } = (N\cap U^{\perp} ) + U =M$$
	where the second to last equality holds because $U\subseteq U^{\perp} $.
\end{proof}

\begin{lemma} \label{cancellation}
	Let $(V,Q)$ and $(W,Q')$ be equivariant quadratic forms such 
	that $(V,Q)$ is  metabolic. 
Then  $(W,Q')$ is metabolic if and only if 
	$(W,Q') \perp (V,Q)$ is metabolic. 
\end{lemma} 

\begin{proof} (see \cite[Lemma 4.1.5]{Annika}) 
	If $(W,Q')$ is metabolic then so is $(W,Q') \perp (V,Q)$. 
	\\
	So assume that $(W,Q') \perp (V,Q)$ is metabolic. 
	Clearly $(V,-Q)$ is metabolic and therefore  
	$X:=(W,Q') \perp (V,Q) \perp (V,-Q)$ is metabolic. 
	Let $$U:= \{ (0,v,v) \mid v\in V\}  \leq X.$$ 
	Then $U$ is isotropic. So by Lemma \ref{incl} 
	there is $M=M^{\perp} \leq X$ 
	with $Q(M) = \{ 0 \}$ such that 
	$$U\leq M = M^{\perp } \leq U^{\perp } = \{ (w,v,v)\mid w\in W, v\in V\}  .$$ 
	Let $\pi : X \to W$ denote the projection onto the first component.
	Then $\pi (M) \subseteq M$ and $\pi(M) $ is
	a self-dual isotropic subspace of $(W,Q')$. 
	In particular $(W,Q')$ is metabolic. 
\end{proof}

\begin{proposition}
	Let $(V,Q)$ and $(W,Q')$ be equivariant quadratic forms. 
	Then 
	$(V,Q)\perp (W,-Q')$ is metabolic if and only if there 
are metabolic modules $(N,H)$, $(M,H')$ such that 
$(V,Q) \perp (N,H) \cong (W,Q') \perp (M,H')$. 
\end{proposition}

\begin{proof} 
	Assume that $(V,Q)\perp (W,-Q')$ is metabolic.
	Then 
	$$(W,Q') \perp ((V,Q)\perp (W,-Q')) \cong (V,Q) \perp ( (W,Q')\perp (W,-Q')) $$ 
	where the two equivariant quadratic forms
	$(V,Q)\perp (W,-Q')$ and $(W,Q')\perp (W,-Q')$ are metabolic. 
	\\
For the opposite direction assume that 
$(V,Q) \perp (N,H) \cong (W,Q') \perp (M,H')$ for metabolic $(N,H)$, $(M,H')$.
Then 
	$$(V,Q) \perp (N,H) \perp (W,-Q') \cong 
	(W,Q') \perp (M,H') \perp (W,-Q') $$
	is metabolic and hence so is $(V,Q)  \perp (W,-Q') $ 
	by Lemma \ref{cancellation}.
\end{proof} 

\subsection{Unique anisotropic representative}

  \begin{lemma} \label{factor}
	  Let $(V,Q)$ be an equivariant quadratic form and 
	  let $U\leq V$ be isotropic.
	  Then 
	  $$\overline{Q}: U^{\perp }/U \to K,\  \overline{Q}(v+U) := Q(v) $$
	  is a well-defined $G$-invariant non-degenerate quadratic form
	  that is Witt equivalent to $(V,Q)$.
  \end{lemma} 

  \begin{proof}
	  Standard computations show that $\overline{Q}$ is well-defined, 
	  $G$-invariant and non-degenerate. 
	  To show that this form is Witt equivalent to $(V,Q)$ we remark 
	  that 
	  $$ \Delta(U^{\perp}) := \{ (v, v+U)\in V\perp U^{\perp}/U \mid v\in U^{\perp} \} $$ is an isotropic subspace of $(V,Q)\perp (U^{\perp}/U,-\overline{Q}) $. 
	  As $\dim(\Delta (U^{\perp})) = \dim(U^{\perp} )$ 
	  and 
	  $$\dim(V) + \dim(U^{\perp} / U)  = \dim(U^{\perp}) + \dim(U^{\vee}) +\dim(U^{\perp}) - \dim(U) = 2\dim (U^{\perp}) $$ 
	  we also get that $\Delta (U^{\perp} ) ^{\perp} = \Delta (U^{\perp} )$.
  \end{proof}

\begin{theorem}  \label{anisoWitt}
	Any class $[(V,Q)] \in \WQ (K,G)$ contains a unique anisotropic 
	representative. 
\end{theorem} 

\begin{proof}
	The existence of an anisotropic representative follows from
	Lemma \ref{factor}. For the uniqueness, let 
	$(V,Q)$ and $(V',Q') \in [(V,Q)]$ be two
	anisotropic modules in the same class of $\WQ(K,G)$. 
	Then $(V,Q) \perp (V',-Q')$ is metabolic, so there is an
	isotropic submodule $U\leq V \oplus V'$ with $U=U^{\perp}$. 
	Clearly $U\meet V$ and $U\meet V'$ are isotropic
	and hence $\{ 0 \}$. So 
	$U=\{ (v,f(v)) \mid v\in V \} $ for some $G$-equivariant 
	isometry $f:(V,Q) \to (V',Q')$. 
\end{proof}

\subsection{Equivariant Witt groups over fields of characteristic not 2} \label{odd}

In this section we briefly recall the 
description of $\WQ(K,G)$  in characteristic not 2. 
Throughout this short section let $K$ be a field of characteristic not 2 
and $G$ be a finite group. 
To obtain an explicit description of the Witt group it suffices 
to enumerate all anisotropic equivariant quadratic forms 
(cf. Theorem \ref{anisoWitt}).
As the characteristic of $K$ is not $2$ the polarization 
$B_Q$ from \eqref{polar} determines the quadratic form $Q$.
In particular the restriction of 
an anisotropic quadratic form to any simple submodule is 
non-degenerate. 
This shows that anisotropic quadratic forms 
are the orthogonal direct sum of simple submodules 
and the Witt group $\WQ(K,G)$ 
can be obtained from 
\cite[Chapter 7]{Scharlau} or \cite[Chapter 4]{Rhiem}.

\begin{lemma}
	Let $(V,Q)$ be an anisotropic equivariant quadratic form.
	Then $(V,Q) = \perp _{j=1}^r (V_j,Q_j) $ for 
	simple $KG$-modules $V_1,\ldots , V_r$.
\end{lemma}

\begin{proof}
	Let $U\leq V$ be a simple submodule of $V$.
Then the restriction $Q_{|U} \neq 0$ because $V$ is anisotropic. 
	Hence also $B_Q(U,U) \neq \{ 0 \}$ so $Q_{|U}$ is non-degenerate 
	and $V=U\perp U^{\perp }$.
	Continue with $U^{\perp }$ instead of $V$ we finally
	achieve an orthogonal decomposition of $V$ into
 equivariant quadratic forms on simple submodules. 
\end{proof}

Now let $V_1,\ldots , V_h$ represent all isomorphism 
classes of simple $KG$-modules admitting 
equivariant quadratic forms $(V_j,Q_j)$. 
Put $D_j:= \End_{KG}(V_j)$. Then $D_j$ is a finite
dimensional $K$-division algebra and the adjoint involution 
of the polarization of $Q_j$ defines an involution $\iota _j$ 
on $D_j$. Denote by $W(D_j,\iota_j)$ the Witt group of 
$\iota _j$-Hermitian forms. 
If $D_j$ is non-commutative then  $\iota _j$ depends on the choice
of $Q_j$ in general. Using equivariant Morita theory
we obtain the following
explicit description of $\WQ(K,G)$:

\begin{theorem} (see \cite[Satz 1.3.8]{habil}, \cite[Section 3.4 (5)]{Scharlaupaper})
	$\WQ(K,G) \cong \bigoplus _{j=1}^h W (D_j,\iota_j) $.
\end{theorem}

\section{Invariants on the equivariant Witt group} \label{invari}

The aim of this section is to define three group homomorphisms 
\begin{center}
$\A : \WQ(K,G) \to \WQ(K) $, $\C:\WQ(K,G) \to \F_2^{{\mathcal S}}$, 
$\D : \ker(\A)  \to \Hom (G,C_2)$ 
\end{center}
on $\WQ(K,G)\cong \ker(\A) \times \WQ(K) $. 
Though these can be defined for general fields, 
they are particularly useful for perfect fields of characteristic 2. 

The {\em orthogonal group} 
of a non-degenerate quadratic space $(V,Q)$
is 
$$O(V,Q):=\{ g\in \GL(V) \mid Q(vg) = Q(v) \mbox{ for all } v\in V \} .$$
The well known {\em Dickson invariant} defines a group homomorphism 
from $O(V,Q) $ to $\{ \pm 1 \} $. If $\Char(K) \neq 2$ then the
Dickson invariant coincides with the determinant. 
In general one defines the Dickson invariant of an element $g\in O(V,Q)$
as $D(g) := (-1) ^{\rk (g-\id_V ) }$.
Then $D : O(V,Q) \to C_2$ is a group homomorphism 
(see \cite[Theorem 11.43]{Taylor}). 

\begin{lemma} (\cite[Lemma 11.58 and Theorem 11.61]{Taylor}) \label{Dickson} 
	Assume that 
	$[(V,Q)] = 0 \in \WQ(K) $ and let $W = W^{\perp} \leq V$ 
	be an isotropic subspace. 
	Then $D(g) = (-1) ^{\dim (W/W\cap Wg) }$  for all $g\in O(V,Q)$.
\end{lemma} 

A second invariant concerns the $KG$-module structure of $V$. 
For this we need the set 
$${\mathcal S} := \{ [S] \mid S \mbox{ is a simple, self-dual } KG\mbox{-module} \} $$ 
 of isomorphism classes of simple, self-dual $KG$-modules.
By the Jordan-H\"older theorem the multiplicity 
$d(V,S)$ of the simple $KG$-module $S$ as a composition factor of $V$
is well defined.

\begin{definition}
	Let $(V,Q)$ be an equivariant quadratic form. 
	\begin{itemize}
		\item[(a)] $\A((V,Q)) := [(V,Q)]\in \WQ(K)$.
		\item[(b)] $\C((V,Q)) := \sum _{[S] \in {\mathcal S}}  
			\overline{d(V,S) } [S] \in \F_2^{{\mathcal S}} $ 
			where 
			$\overline{\phantom{x}} $ denotes the reduction 
			modulo $2$. 
		\item[(c)] $\D((V,Q)): (g \mapsto D(g_V)) \in \Hom (G,C_2), $  where $g_V \in O(V,Q) $ 
			is the endomorphism of $V$ describing
the action of $g\in G$. 
	\end{itemize}
\end{definition}

\begin{theorem} \label{ACD}
	The maps
	$\A$ and $\C$ are well defined group homomorphisms on $\WQ  (K,G)$. 
	\\
	The map $\D$ is a well defined group homomorphism on $\ker(\A)$.
\end{theorem} 

\begin{proof}
	Clearly the
	forgetful homomorphism $\A:\WQ(K,G) \to \WQ(K)$ 
	is a well defined group homomorphism. 
	\\
	To see that $\C$ is well defined on $\WQ(K,G)$ it is enough to remark that for an isotropic submodule $U\leq V$ (as in Lemma \ref{factor}) 
	the module $V/U^{\perp} \cong U^{\vee }$. 
	So any self-dual composition factor of $U$ is also a composition 
	factor of $V/U^{\perp }$ and hence it appears with odd multiplicity
	in $V$ if and only if it appears with odd multiplicity in $U^{\perp }/U$. 
	Clearly $\C $ is compatible with the addition on $\WQ(K,G)$
	defined by orthogonal direct sums. 
	\\
	For the Dickson invariant we use the definition of $D$ 
	from Lemma \ref{Dickson} as 
	$D(g) = (-1) ^{\dim(W/(W\cap Wg)) } $ for any 
	isotropic subspace
	$W=W^{\perp } \leq (V,Q)$. 
	If 
	$ U \leq U^{\perp } \leq V$ is as in Lemma \ref{factor} 
	then $U$ is contained in a maximal isotropic subspace $W$ of $V$ 
	and $W/U \leq U^{\perp } / U$ is a maximal isotropic subspace
	of $(U^{\perp}/U ,\overline{Q})$.
	As $U$ is $G$-invariant we have 
	$$\dim (W/(W\cap Wg)) = \dim ( (W/U) / (W/U \cap (W/U)g)) \mbox{ for all } g\in G$$ 
	and hence $\D $ is also well defined. 
	Again the compatibility of $\D $ with orthogonal direct sums 
	is clear.
\end{proof}

\begin{remark} \label{split} 
	The group homomorphism $\WQ(K) \to \WQ(K,G)$, 
	$[(V,Q)] \mapsto [(V,Q)]$, where $G$ acts trivially on $V$,
	is injective and establishes a decomposition  
	$\WQ(K,G) = \ker(\A) \times \WQ(K) $.
\end{remark}

\section{Anisotropic equivariant quadratic forms} 

In this section let  $K$ be a finite field of characteristic 2. 
Then $\wp(K):=\{ a^2+a \mid a\in K \} $ is a subgroup of the 
additive group 
$$(K,+) = \wp (K) \cup \alpha + \wp(K) $$ 
where $\alpha \in K$ is any element for which the polynomial $X^2+X+\alpha \in K[X]$ 
is irreducible. The Witt group $\WQ(K)$ of quadratic forms over $K$ 
consists of two classes, the trivial one and $[N(K)]$, 
where $N(K) = \langle f,e \rangle $ with $Q(f)=\alpha $, $Q(e) =1 $, 
$B_Q(e,f) = 1$ is the norm form of the quadratic extension of $K$.
This is the unique non-zero anisotropic quadratic form over $K$ 
(cf. \cite[Section 12]{Kneser}).

We also note that squaring is a field automorphism of $K$ thus
every element of $K$ has a square root. 
For a $K$-space $V$ we denote by $V^{(2)}$ the $K$-space with 
the same underlying abelian group $V$ where $K$ acts by 
$K\times V^{(2)} \to V^{(2)}, (a,v) \mapsto a^2 v $. 

Let $G$ be a finite group.
Denote by $T$ the trivial $KG$-module, i.e. $T=K$ and $vg = v $ for all $v\in T,
g\in G$. 
Then squaring yields a $KG$-module isomorphism $T \cong T^{(2)} $.

If a simple $KG$-module $V$ admits a $G$-invariant non-degenerate
quadratic form, then $V$ is called of {\em orthogonal type}. 
Of course orthogonal $KG$-modules are self-dual. 
If a non-trivial self-dual $KG$-module $W$ does not admit a 
non-zero quadratic form then $W$ is called of {\em symplectic type}. 
The set ${\mathcal S}$ of self-dual simple $KG$-modules 
from Section \ref{ACD} is hence of
the form ${\mathcal S} = {\mathcal S}_0 \disj \{ [T ] \} $ where 
$${\mathcal S}_0 = \{ [V] \mid V \mbox{ self-dual, simple, orthogonal} \} 
\disj \{ [W] \mid W \mbox{ self-dual, simple, symplectic} \} .$$

\begin{theorem} (see also \cite[Section 4.1]{Annika}) \label{simple} 
	Let $V$ be a simple self-dual $KG$-module. 
	\begin{itemize} 
		\item[(a)]
	There is a non-zero $G$-invariant symmetric 
	bilinear form on $V$.
\item[(b)]
	Any non-zero $G$-invariant 
	symmetric bilinear form $B$ on $V$ is non-degenerate. 
\item[(c)] 
	If $V$ is not the trivial $KG$-module then 
	any $G$-invariant symmetric bilinear form $B$ on $V$ is symplectic,
	i.e. $B(v,v) = 0$ for all $v\in V$.
\item[(d)] Any two non-zero $G$-invariant symmetric bilinear
	forms on $V$ are $KG$-isometric. 
\item[(e)] 
If $V$
admits a non-zero $G$-invariant quadratic form $Q$ then 
	either $V\cong T$ and $B_Q=0$ or $Q$ is non-degenerate.
	In both cases the non-zero $G$-invariant quadratic form is
unique up to $KG$-isometry. 
\end{itemize} 
\end{theorem}

\begin{proof}
	(b)
Let $B:V\times V \to K$ be a non-zero symmetric 	
$G$-invariant form. Then $V^{\perp} \leq V$. As $V$ is simple
we have $V^{\perp }=\{0\}$ and hence $B$ is non-degenerate. 
\\
	(a) 
	Let $f: V\to V^{\vee }$ be a $KG$-module isomorphism.
	Define $\beta :V\times V \to K $ by $\beta (v,w):=f(w)(v) $. 
	Then either $\beta (v,w) = \beta (w,v) $ for all $v,w\in V$ 
	and $\beta $ is a non-degenerate symmetric $G$-invariant 
	bilinear form or 
	$B:V\times V \to K, B(v,w):=\beta(v,w)+ \beta (w,v)$ is a 
	non-zero symmetric bilinear form on $V$.
	As $V$ is simple $B$ is non-degenerate by (b).
\\
	(c) The map $Q_B:V\to T^{(2)}$, $v\mapsto B(v,v) $ is 
	a $KG$-module homomorphism, because
	$$Q_B(v+aw) = B(v+aw,v+aw) = B(v,v)+a^2B(w,w) = Q_B(v) + a^2Q_B(w) $$ 
	for all $v,w\in V, a\in K$. 
	As $V\not\cong T$ and $V$ is simple $Q_B=0$. 
	\\
	(d)
	As $V$ is simple, its endomorphism ring 
	$E:=\End_{KG}(V)$ is again a finite field of characteristic 2.
	Moreover the adjoint involution $\ ^{ad} $ of $B$
	defined by 
	$$B(ve,w) = B(v,we^{ad}) \mbox{ for all } e\in E, v,w \in V $$ 
	defines a $K$-linear field automorphism of $E$ of order 1 or 2.
	The space $E^+:=\{ e\in E \mid e = e^{ad} \}$ is a subfield of $E$ and
the map 
	$E\to E^+, e\mapsto e e^{ad} $ is either the norm or squaring, 
	in particular it is surjective. 

All non-degenerate $G$-invariant symmetric bilinear forms on $V$ are of the 
	form 
	$$sB:V\times V\to K, sB (v,w):=B(v,ws) \mbox{ 
	for some non-zero } s\in E^+.$$ There is $e\in E$ such that 
	$s = e e^{ad}$, so 
	$sB(v,w) = B(ve,we)$ and multiplication by $e$ defines a 
	$KG$-isometry between $(V,B)$ and $(V,sB)$. 
\\
	(e) Let $Q$ be a non-zero $G$-invariant quadratic form on $V$.
	If $B_Q = 0$ then $Q:V\to T^{(2)}$ is linear and hence 
	$V\cong T$. 
	As squaring is surjective, the
	 non-zero quadratic form on $T$ is unique.
	\\
        Now 
	assume that $V\not \cong T$. 
	Then $B_Q$ is non-zero and therefore  non-degenerate by (b).  
	\\
	Let $Q'$ be a second $G$-invariant quadratic form on $V$. 
	By (d) we may assume that $B_Q=B_{Q'}$. 
	So $Q-Q' : V \to T^{(2)}$ is linear.  As $V\not \cong T$ we have
	$Q=Q'$. 
\end{proof}

\begin{corollary}\label{mult}
	Let $V$ be a simple $KG$-module and $Q,Q'$ be two non-degenerate
	$G$-invariant quadratic forms on $V$.
	Then $[(V,Q)] + [(V,Q')] = 0 \in \WQ(K,G)$. 
\end{corollary}

\begin{proof}
	Let $f:(V,Q) \to (V,Q')$ be a $KG$-isometry (see Theorem \ref{simple} (e)). Then $U:=\{ (v,f(v) ) \mid v\in V \} \leq (V,Q) \perp (V,Q') $ 
	is an isotropic $KG$-submodule with $U=U^{\perp }$.
\end{proof}

As the bilinear forms on the simple modules are unique up to $KG$-isometry
and the endomorphism ring of the direct sum of distinct simple $KG$-modules
is the direct sum of the endomorphism rings of the simple summands
Theorem \ref{simple} (d) also implies the following corollary. 

\begin{corollary} \label{unique}
	If $W$ is a direct sum of pairwise non-isomorphic simple 
	self-dual $KG$-modules then there is a non-degenerate 
	symmetric $G$-invariant bilinear form  $B:W\times W \to K$.
	Such a form is unique up to $KG$-isometry.
\end{corollary}

\begin{theorem}\label{aniso}
Let $(V,Q)$ be an anisotropic $KG$-module.  \\
Then the socle of $V$ is  
${\mathcal T} \perp V_0 $, 
where $V_0$ is the orthogonal sum of 
pairwise non isomorphic simple $KG$-modules of 
orthogonal type and  either
\begin{itemize}
\item[(i)] ${\mathcal T} = \{ 0 \} $ and $V=V_0$.
\item[(ii)] ${\mathcal T} \cong T\oplus T $ 
and $Q_{|{\mathcal T}}$ is the 
unique anisotropic $2$-dimensional quadratic form $N(K)$ over $K$. 
Then $V = V_0 \perp {\mathcal T}$. 
\item[(iii)] ${\mathcal T} = T = \langle e \rangle $ with $Q(e) = 1$ 
	and $V=V_0 \perp R$ for some indecomposable $KG$-module $R$ with 
	socle $T$.
\end{itemize}
\end{theorem}

\begin{proof}
   Let $U\leq V$ be a simple submodule. Then $Q(U) \neq \{ 0 \}$ and
	hence by part (e) of Theorem \ref{simple} either
        the restriction of $Q$ to $U$ is non-degenerate or
 $U=T$. 
        In the first case $V$ splits as $V=U\perp U^{\perp }$. 
	Continuing like this, we arrive at a decomposition 
	$V=V_0 \perp V_0^{\perp }$ where $V_0$ is an orthogonal
	sum of simple orthogonal $KG$-modules. 
As $V_0$ is anisotropic, it is multiplicity free by Corollary \ref{mult}.
\\
	Replacing $V$ by $V_0^{\perp}$ we hence may assume that the
	socle ${\mathcal T}$ of $V$ is the direct sum of trivial $KG$-modules. 
	If ${\mathcal T}=\{ 0\}$ then we are in case (i). 
	\\
	Now ${\mathcal T}$ is an anisotropic quadratic $K$-space, 
	so by \cite[Section 12]{Kneser} 
	 $\dim _K({\mathcal T}) \leq 2$ and if 
	$\dim_K({\mathcal T}) = 2$ then 
	${\mathcal T} \cong N(K)$.
	In particular ${\mathcal T} $ 
	is non-degenerate and hence splits as an orthogonal summand.
	So then $V_0^{\perp } = {\mathcal T}$. 
	\\
	If $\dim_K({\mathcal T}) = 1$, then ${\mathcal T} = T = \langle e \rangle $ with $Q(e) = 1$ and $R=V_0^{\perp }$ is indecomposable as it has
	a simple socle.
\end{proof}

We now analyse the module $R$ from Theorem \ref{aniso} (iii). 
Related ideas can be found in \cite[Theorem 1.3]{GowWillems}.

\begin{theorem}\label{anisoR}
	Let $R$ and $e$ be as in Theorem \ref{aniso} (iii)
	and put $W:= \langle e \rangle ^{\perp } / \langle e \rangle   $ 
	where $\langle e \rangle ^{\perp }$ is computed in $R$.
	Then either $W=\{ 0\}$ or 
	$W \cong \oplus _{j=1}^r W_j $ 
	is the direct sum of  pairwise non-isometric simple 
	$KG$-modules $W_j$ of symplectic type. 
\end{theorem} 

\begin{proof}
	\begin{itemize}
		\item $W$ is the orthogonal sum of simple $KG$-modules: \\
	As $B_Q(e,e)=0$, the bilinear form $B_Q$ 
	defines a non-degenerate bilinear form $B$ on $W$.
	 Let $U$ be a simple submodule of $W$ and $\tilde{U} \leq \langle e\rangle ^{\perp} $
        denote its full preimage in $\langle e\rangle ^{\perp }$.
	If the restriction $B_{|U}$ of $B$ to $U$ is zero, then
	$Q:\tilde{U} \to T^{(2)}$ is a $KG$-module homomorphism.
 As $Q(e) \neq 0 $ we get
        $\tilde{U} = \ker(Q) \perp \langle e \rangle $.
        This contradicts the fact that the socle of $R$ is $\langle e \rangle$.
	So $B_{|U}$ is non-degenerate and
        $W = U \perp U^{\perp }$.
		\item 
        We now show that $U$ is symplectic. 
			\\
			Otherwise there is a
        $G$-invariant non-degenerate quadratic form $F:U\to K$
        such that $B = B_F$. Extend $F$ to a quadratic form on $\tilde{U}$
	with $F(e) = 0 $. Then $B_{Q+F} = 0$ on $\tilde{U}$ and hence
	$Q+F :\tilde{U} \to T^{(2)}$ is a $KG$-homomorphism giving the
        same contradiction as before.
\item $W$ is multiplicity free. \\
        Assume that there is a submodule $U'\leq U^{\perp }$ that
        is isomorphic to $U$ and choose an isometry
	$f:(U,B_{|U}) \to (U',B_{|U'}) $ (see Theorem \ref{simple} (d)). 
	Then $B$ is identically zero on 
        $U'':=\{ u+f(u) \mid u \in U \} \leq W$. Clearly $U\cong U''$ and
	as before $Q: \tilde{U''} \to T^{(2)}$ is an epimorphism with 
        $\tilde{U''} = \ker(Q) \perp \langle e \rangle $. 
	\end{itemize} 
\end{proof}
	
\begin{theorem}\label{anisoR0}
	If $W=\{ 0\} $ in Theorem \ref{anisoR}
	then $R$ has a $K$-basis $(f,e)$ with $B_Q(f,e) = 1$ and 
	either $Q(f) = 0$ or $Q(f) = \alpha \not\in \wp(K) $.
With respect to this basis 
	 $G$ acts on $R$ as 
$\langle \left( \begin{array}{cc} 1 & 1 \\ 
	0 & 1 \end{array} \right) \rangle $. 
	For each epimorphism $\tau $ of $G$ onto a group of order $2$ 
	there are two such modules $R$, $R^+(\tau )$ and $R^-(\tau )$,
	where the underlying $2$-dimensional quadratic space of 
	$R^+(\tau )$ is the hyperbolic plane  (i.e. $Q(f) = 0$)  
	and for $R^-(\tau )$ this is the 
unique anisotropic $K$-space $N(K)$ (i.e. $Q(f) = \alpha $).
\end{theorem}

\begin{proof} As $Q$ is non-degenerate and $Q(e) =1$ the module 
 $R$ has a basis 
	$(f,e)$  with $B_Q(e,f) =1$.
	These two conditions uniquely
	determine $e\in \soc(R)$ and the class $f +\langle e \rangle $.
	We have $Q(f+ae) = Q(f) + (a+a^2) $ so we can achieve that 
	$Q(f) \in \{ 0,\alpha \}$ where $X^2+X+\alpha \in K[X] $
	is irreducible. As $g\in G$ fixes $e$ it also fixes the class  
	$f+\langle e \rangle $ and hence either fixes $f$ or maps 
	$f$ to $f+e$. 
\end{proof}

 \section{The quadratic envelope of a symplectic $KG$-module} 

 We keep the assumption that $K$ is a finite field of characteristic 
 $2$.
Our considerations are inspired by 
 \cite[Theorem 11.9]{Taylor} that establishes 
 an isomorphism between $O_{2m+1}(K) $ and $ \Sp_{2m}(K)$.  
 In our context the following lemma seems to be easier to use:
 
 \begin{lemma}\label{lem:orthogonalisation}
	 Let $(R,Q)$ be a non-degenerate 
	 quadratic space of dimension $2m+2$ over $K$ of maximal Witt index $m+1$. Let $e\in R$ be such that $Q(e)=1$. 
	 Then 
	 $$S(e):= \{ g \in O(R,Q) \mid eg =e , D(g) = 1 \} \cong \Sp _{2m}(K) .$$
 \end{lemma}

 \begin{proof}
Let $$(f,w_1,\ldots , w_m,v_1,\ldots , v_m,e )$$ be a basis of $R$, 
	 such that $\langle e,f \rangle $, $\langle v_i,w_i \rangle $ 
	 ($1\leq i \leq m $) 
	 are pairwise orthogonal hyperbolic planes, $Q(f)=0$, $Q(e) =1$, $B_Q(e,f) = 1$, 
	 $B_Q(v_i,w_i) =1 $, and 
	 $Q(v_i) = Q(w_i) = 0$ for all $i=1,\ldots ,m$. 
	 Any element $g\in O(R,Q)$ with $eg=e$ also stabilises $\langle e\rangle ^{\perp} $ 
	 and the class $f+\langle e \rangle ^{\perp }$ and hence its matrix is of the form
	$$g = 
	\left( \begin{array}{cccc} 
		1 & a & b & x \\ 
		0 & A & B & c \\
		0 & C & D & d \\
	0 & 0 & 0 & 1 \end{array} \right) .$$ 
	We hence obtain a group homomorphism 
	$$\varphi : S(e) \to \Sp_{2m}(K): g \mapsto 
\left( \begin{array}{cc} A & B \\ C & D \end{array}\right) .$$
	 As $Q(w_i g) = 0$ we have $c_i^2 = (AB^{tr})_{ii}$. 
	 Similarly $Q(v_ig)=0 $ implies that $d_i^2=(CD^{tr})_{ii} $. 
	 So the restriction $g'$ of $g$ to  $\langle e \rangle ^{\perp }$ is uniquely determined 
	 by $\varphi(g) $.
	 Now $g': \langle e \rangle ^{\perp } \to \langle e \rangle ^{\perp } $ is an isometry so by Witt's extension
	 theorem there is $g\in O(R,Q)$ such that $g_{|\langle e \rangle^{\perp}} = g'$. 
	 The conditions that $B_Q(fg,w_ig) =0$
	 and  $B_Q(fg,v_ig)=0$ for all $i$ yield 
	 $$\left(\begin{array}{c} c \\ d \end{array}\right) = 
		 \left( \begin{array}{cc} A & B \\ C & D \end{array}\right) 
			 \left(\begin{array}{c} b^{tr} \\ a^{tr} \end{array}\right).$$ 
	 and hence uniquely determine $a,b \in K^m$.
Now Witt's extension theorem implies that 
	 $0=Q(fg) = ab^{tr} + x^2+x $ has a solution $x\in K$, so 
	 $ab^{tr} \in \wp (K)$. 
	 In fact the equation $ab^{tr} + x^2+x =0$ then has two solutions,
	 say $x_0$ and $x_0+1$.
	 So 
	$$g =  g_0 := 
	\left( \begin{array}{cccc} 
		1 & a & b & x_0 \\ 
		0 & A & B & c \\
		0 & C & D & d  \\
	0 & 0 & 0 & 1 \end{array} \right) \mbox{ or } 
	 g = g_1 = g_0 h \mbox{ with }  h = 
	\left( \begin{array}{cccc} 
		1 & 0 & 0 & 1 \\ 
		0 & I_m & 0 & 0 \\
		0 & 0 & I_m  & 0  \\
	0 & 0 & 0 & 1 \end{array} \right) $$ 
	 where $a,b,c,d\in K^m$ are uniquely determined by $\varphi (g)$.
		 The Dickson invariant $D(h) = -1$ so exactly 
		 one of $g_0$ or $g_1$ has trivial Dickson invariant. 
		 Therefore $\varphi $ is the desired isomorpism.
 \end{proof}

Now let $G$ be a finite group and let
$W_1,\ldots , W_r$ be pairwise non-isomorphic simple symplectic $KG$-modules 
 and put $W:=W_1\oplus \ldots \oplus W_r $.  
 We assume that $W \neq \{ 0 \}$.

By  Corollary \ref{unique} there is a unique non-degenerate 
$G$-invariant symplectic bilinear form $B$ on $W$. 
Then the action of $G$ on $W$ yields a homomorphism 
$$\rho_W: G \to \Sp(W) \cong \Sp_{2m}(K) = S(e) \leq O(R,Q)$$ 
with $(R,Q)$ as in Lemma \ref{lem:orthogonalisation}.

\begin{definition} \label{orthogonalisation}
	The equivariant quadratic form $(R,Q)$ is called 
	the {\em quadratic envelope} 
	$R^{+}(W)$ of the symplectic $KG$-module $W$.
\end{definition}

We summarize the properties of $R^+(W)$ in the following proposition:

\begin{proposition} \label{simplerump}
	$R^{+}(W) = (R,Q) $ is an anisotropic equivariant quadratic form. 
 \begin{itemize}
	 \item[(a)] $\soc (R^+(W)) = \langle e \rangle \cong T$ with 
		 $Q(e) = 1$. 
	 \item[(b)] $\langle e \rangle^{\perp}/\langle e \rangle \cong W$. 
	 \item[(c)] $\A (R^+(W)) = 0 $. 
	 \item[(d)] $\D (R^+(W)) = 1 $.
	 \item[(e)] $\C (R^+(W)) = \sum_{j=1}^r [W_j] $.
 \end{itemize} 
\end{proposition}

\begin{proof}
	(a) By construction $\langle e \rangle \cong T $ is a $KG$-submodule of 
	the socle of $R$. 
	Assume first that a direct summand of $W$ is a summand $W_j$ 
	of $\soc (R)$. 
	As $W_j$ is a symplectic $KG$-module, the restriction of 
	$Q$ to $W_j$ is 0. Now $W_j$ is self-dual and occurs with 
	multiplicity 1 in $R^+(W)$, so this implies that 
	$W_j$ is in the radical of $Q$, a contradiction. 
	If $\soc(R) \cong T\oplus T$ then the restriction of 
	$Q$ to $\soc(R)$ is non-degenerate and $\soc(R)$ 
	splits as an orthogonal direct summand, implying that 
	$R = T\oplus T$ and hence $W=\{0\}$, contradicting our assumption.
	\\
	(b), (c), and (e) are clear by construction and (d) follows from 
	the choice of $g = g_0$ or $g_0 h$ in the proof of Lemma \ref{lem:orthogonalisation} to guarantee that the Dickson invariant of $g$ be trivial. 
\end{proof}

The construction of the quadratic envelope shows that every 
simple $KG$-module of symplectic type has a non-trivial 
extension with the trivial module. 
The following important consequence is well known.

\begin{corollary}  (cf. \cite[Proposition 2.4]{SinWillems})
	All simple self-dual $KG$-modules of symplectic type lie in the 
	principal block. 
\end{corollary}

\section{The Witt group of $KG$} 

We now use the invariants of the Witt group defined in Section \ref{invari} 
to describe the Witt group of equivariant quadratic forms 
for a finite group $G$ and a finite field $K$ of characteristic 2. 
Recall that the Witt group of quadratic forms 
$\WQ(K) = \{ 0, [N(K)] \} $ is a group of order $2$ (see \cite[Satz 12.4]{Kneser}). 
Recall the definition of 
${\mathcal S}_0 := {\mathcal S} \setminus \{ [T] \}  $ 
as the set of non-trivial self-dual simple $KG$-modules.

\begin{theorem}
	$\WQ(K,G) \cong \F_2^{{\mathcal S}_0} \times \WQ(K) \times \Hom (G,C_2)  
	\cong C_2^{s+t} $, where $s= |{\mathcal S} |$ and 
	$t$ is the $2$-rank of $G/G'$.
\end{theorem} 

\begin{proof}
	\begin{itemize}
		\item[(a)] By Remark \ref{split} we have 
			$$\WQ(K,G) = \ker(\A) \times \WQ (K) = \ker(\A) \times 
			\langle [N(K)] \rangle .$$
			Clearly $[N(K)] \in \ker(\C )$.
		\item[(b)] We now show that 
	$\C :\WQ(K,G) \to \F_2^{{\mathcal S}_0}$ is surjective:
			\\
	Any 
orthogonal simple $KG$-module  
	$(V,Q)$ is an anisotropic representative of its 
	class $[(V,Q)] \in \WQ(K,G) $. 
	We have $\C ([(V,Q)] ) = [V]$. 
\\	
	For any symplectic simple $KG$-module $W$
	 Proposition \ref{simplerump} constructs an anisotropic 
	 equivariant quadratic form $R^+(W)$ with 
	 $\C([R^+(W)] ) = [W] $. 
 \item[(c)] As $[N(K)]\in \ker(\C) $ we now conclude
	 that $\A \times \C : \WQ (K,G) \to \WQ(K) \times \F_2^{{\mathcal S}_0} $ is surjective and split. The 
		subgroup 
			$$\langle [(V,Q)],[R^{+}(W)],[N(K)] \mid V \mbox{ simple orthogonal, } W \mbox{ simple, symplectic} \rangle \cong C_2^{s} $$
			generates a complement of 
			$\ker(\A ) \meet \ker(\C) $.
\item[(d)]
	Theorem \ref{anisoR0}  shows that 
			$$\ker(\A ) \meet \ker(\C) =\langle [R^+(\tau _j )] \mid  j=1,\ldots ,t \rangle \cong \Hom (G,C_2) .$$
	\end{itemize} 
\end{proof}

%



\begin{thebibliography}{99}

	\bibitem{Dress}  Andreas Dress, {Induction and structure theorems 
		for orthogonal representations of finite groups}, 
		Annals of Mathematics,
		Second Series,  {\bf 102} (1975) 291--325.
	\bibitem{GowWillems} 
Roderick Gow, Wolfgang Willems, 
Methods to decide if simple self-dual modules over fields of characteristic 2
are of quadratic type.
		J. Algebra {\bf 175} (1995) 1067--1081.

\bibitem{Kneser} Martin Kneser, {\em Quadratische Formen.} 
Neu bearbeitet und herausgegeben in Zusammenarbeit mit Rudolf Scharlau.
		Springer (Berlin) (2002)
	\bibitem{Annika} Annika Meyer, 
		{\em Automorphism groups of self-dual codes.}
Dissertation, RWTH Aachen University (2009)
\bibitem{habil}  Gabriele Nebe, {\em Orthogonale Darstellungen
                endlicher Gruppen und Gruppenringe}, Habilitationsschrift,
                Aachener Beitr\"age zur Mathematik  26 (1999) Verlag Mainz, Aachen.
	\bibitem{Scharlaupaper} 
		Heinz-Georg Quebbemann, Winfried Scharlau, M. Schulte, 
		{Quadratic and Hermitian Forms in Additive and Abelian Categories.} J. Algebra {\bf 59} (1979) 264--289.
\bibitem{Rhiem} 
Carl Rhiem, {\em Introduction to Orthogonal, Symplectic and
	Unitary Representations of Finite Groups} 
	AMS Fields Institute Monographs (2011)
\bibitem{Scharlau} Winfried Scharlau, {\em Quadratic and Hermitian forms.}
	Grundlehren der Mathematischen Wissenschaften
		{\bf 270}. Berlin etc.: Springer-Verlag (1985).
		\bibitem{SinWillems}
                Peter Sin, Wolfgang Willems,
                $G$-invariant quadratic forms,
                J. Reine Angew. Math. {\bf 420} (1991) 45--59.
	\bibitem{Taylor} 
		Don Taylor, {\em The geometry of the classical groups.} 
		Heldermann Verlag (1992)
\end{thebibliography}
\end{document}